\documentstyle[fleqn,12pt]{article}
\begin{document}
\pagenumbering{arabic}\setcounter{page}{1}
\pagestyle{plain}\baselineskip=18pt

\begin{center}
{\bf DIFFERENTIAL GEOMETRY OF THE LIE SUPERALGEBRA OF THE QUANTUM
SUPERPLANE}
\end{center}

\vspace{1cm}
\begin{center} Salih \c Celik \footnote{E-mail address: sacelik@yildiz.edu.tr}

Yildiz Technical University, Department of Mathematics, \\
34210 Davutpasa-Esenler, Istanbul, TURKEY. \end{center}

\vspace{2cm}
{\footnotesize Differential geometry of the quantum Lie superalgebra of the
extended quantum superplane and its Z$_2$-graded Hopf algebra structure is
obtained. Its Z$_2$-graded dual Hopf algebra is also given.}

\noindent
{\bf 1. Introduction}

\noindent Noncommutative differential geometry has attracted
considerable interest both mathematically  and also from
theoretical physics side over the past decade. Especially, there
is much activity in differential geometry on quantum groups. For
references to the literature for quantum groups we refer to the
recent book by Majid [1]. The basic structure giving a direction
to the noncommutative geometry is a differential calculus on an
associative algebra. A noncommutative differential calculus on
quantum groups has been introduced by Woronowicz [2]. Wess and
Zumino [3] has been reformulated to fit this general theory, in
less abstract way. Some other methods to define a differential
geometric structure (or a De Rham complex) on a given
noncommutative associative algebra or to construct a
noncommutative geometry on a quantum group have been proposed and
investigated by several authors [4-9].

It is known that, in order to construct a noncommutative differential calculus
on quantum groups and Hopf algebras, one takes into consideration the
associative algebra of functions on the group. The starting point of this work
is its Lie algebra. We present here a differential calculus on the Lie
superalgebra of the associative algebra of functions on the extended
$q$-superplane.

The paper is organized as follows. In the second section, we state
briefly the properties of Z$_2$-graded quantum superplane which
are described in Ref. 10. In the third section, using the fact
that an element of a Lie group can be represented by exponential
of an element of its Lie algebra, we shall write the generators of
the extended $q$-superplane as exponential of some elements. The
obtained new elements which are to be generators of the Lie
superalgebra [11]. We give a differential calculus on the Lie
superalgebra and its Hopf algebra structure. We also obtain its
Z$_2$-graded dual Hopf algebra.

\noindent
{\bf 2. Review of hopf algebra ${\cal A}$}

\noindent The quantum superplane [12] is defined as an
associative algebra whose the even coordinate $x$ and the odd
coordinate $\theta$ satisfying
$$ x \theta - q \theta x = 0 \qquad \theta^2 = 0 \eqno(1) $$
where $q$ is a nonzero complex deformation parameter. This algebra is known as
the algebra of polynomials over the quantum superplane and we shall denote by
${\cal A} = Fun_q(R(1|1))$.

We know that the algebra ${\cal A}$ is a Z$_2$-graded Hopf algebra
with the following costructures [10]: the coproduct $\Delta: {\cal
A} \longrightarrow {\cal A} \otimes {\cal A}$ is defined by
$$\Delta(x) = x \otimes x, \qquad
  \Delta(\theta) = \theta \otimes x + x \otimes \theta. \eqno(2a)$$
The counit $\epsilon: {\cal A} \longrightarrow {\cal C}$ is given
by
$$\epsilon(x) = 1, \qquad \epsilon(\theta) = 0. \eqno(2b)$$
We extend the algebra ${\cal A}$ by including the inverse of $x$
which obeys
$$x x^{-1} = 1 = x^{-1} x.$$
If we extend the algebra ${\cal A}$ by adding the inverse of $x$
then the algebra ${\cal A}$ admits a coinverse $\kappa: {\cal A}
\longrightarrow {\cal A}$ defined by
$$\kappa(x) = x^{-1}, \qquad \kappa(\theta) = - x^{-1} \theta x^{-1}. \eqno(2c)$$
It is not difficult to verify the following properties of
costructures:
$$(\Delta \otimes \mbox{id}) \circ \Delta =
  (\mbox{id} \otimes \Delta) \circ \Delta, $$
$$m_{\cal A} \circ (\epsilon \otimes \mbox{id}) \circ \Delta
  = m_{\cal A} \circ (\mbox{id} \otimes \epsilon) \circ \Delta, \eqno(3)$$
$$m_{\cal A} \circ (\kappa \otimes \mbox{id}) \circ \Delta = \epsilon
  = m_{\cal A} \circ (\mbox{id} \otimes \kappa) \circ \Delta $$
where id denotes the identity mapping and $m_{\cal A}$ is the multiplication
map $$m_{\cal A} : {\cal A} \otimes {\cal A} \longrightarrow {\cal A}, \qquad
  m_{\cal A}(a \otimes b) = ab. $$
The multiplication in ${\cal A} \otimes {\cal A}$ follows the rule
$$(A \otimes B) (C \otimes D) = (-1)^{\hat{B} \hat{C}} AC \otimes BD,
  \eqno(4)$$
where $\hat{f}$ denotes the Z$_2$-grading of $f$.

\noindent
{\bf 3. Differential calculus on the Lie superalgebra}

\noindent
It is known that an element of a Lie group can be represented by exponential
of an element of its Lie algebra. Using this fact, one can define the
generators of $\cal A$ as
$$x = e^u, \qquad \theta = e^u \eta. \eqno(5)$$
Let
$$ q = e^h. \eqno(6)$$
Then we obtain the relations
$$[u,\eta] = h \eta, \qquad \eta^2 = 0, \eqno(7)$$
where
$$[a,b]_\pm = a b \pm ba. $$
These are the relations of a Lie superalgebra and we shall denote it by
$\cal L(A)$. The Z$_2$-graded Hopf algebra structure of $\cal L(A)$ can be
read off from (2),
$$\Delta(u) = u \otimes 1 + 1 \otimes u, \qquad
  \Delta(\eta) = \eta \otimes 1 + 1 \otimes \eta, $$
$$\epsilon(u) = 0, \qquad \epsilon(\eta) = 0, \eqno(8)$$
$$\kappa(u) = - u, \qquad \kappa(\eta) = - \eta. $$
An interesting case is that these costructures are the same with the Hopf
algebra structure of one-forms on $\cal A$ which is given in Ref. 10.

We want build up a noncommutative differential calculus on the Lie
superalgebra $\cal L(A)$. This may be involve functions on the Lie
superalgebra $\cal L(A)$, differentials and differential forms. So we
have to define a linear operator {\sf d} which acts on the functions of
the elements of $\cal L(A)$.

In order to establish a noncommutative differential calculus on the algebra
$\cal L(A)$, we assume that the commutation relations between the
elements of $\cal L(A)$ and their differentials are of the following form:
$$u ~{\sf d}u = A_{11} {\sf d}u ~u + B_1 {\sf d} u + B_2 {\sf d} \eta, $$
$$u ~{\sf d}\eta = A_{12} {\sf d}\eta ~u + B_3 {\sf d}u + B_4 {\sf d}\eta,$$
$$\eta ~{\sf d}u = A_{21} {\sf d}u ~\eta + B_5 {\sf d}u + B_6 {\sf d} \eta,$$
$$\eta ~{\sf d}\eta = A_{22} {\sf d}\eta ~\eta + B_7 {\sf d}u +
  B_8 {\sf d} \eta. \eqno(9)$$
The coefficients $A_{ij}$ and $B_i$ will be determined in terms of the
"new" deformation parameter $h$. To find them we shall use the
consistency of calculus.

We first note that the following properties of the exterior differential
{\sf d}: the nilpotency
$${\sf d}^2 = 0, \eqno(10\mbox{a})$$
and the Z$_2$-graded Leibniz rule
$${\sf d}(f g) = ({\sf d} f) g +(-1)^{\hat{f}} f ({\sf d} g). \eqno(10\mbox{b})$$
From the consistency conditions
$${\sf d}(u \eta - \eta u - h \eta) = 0, \qquad {\sf d}(\eta^2) = 0 $$
we find
$$A_{12} = 1, \qquad B_3 + B_5 = 0, \qquad A_{22} = 1, \eqno(11\mbox{a})$$
$$A_{21} = -1, \qquad B_4 + B_6 = h, \qquad B_7 = 0 = B_8. $$

Similarly, from
$$(u \eta - \eta u - h \eta) {\sf d}u = 0, \qquad
  (u \eta - \eta u - h \eta) {\sf d}\eta = 0 $$
one has
$$B_2 = 0 = B_3, \qquad A_{11} B_5 = 0, \qquad
  (B_1 - A_{11}) B_5 = 0 = (1 - A_{11})B_6, \eqno(11\mbox{b})$$
$$(B_1 - B_4 + h)B_6 = 0, \qquad (B_1 +h)B_5 = B_3 B_6. $$
The system (11) has many solutions and we shall only discuss one of them
below. Most of the coefficients in the relations (9) are already determined.
The remaining coefficients can be determined from the following equations
$$ B_4 + B_6 = h, \qquad (1 - A_{11})B_6 = 0, \qquad
  (B_1 - B_4 + h)B_6 = 0. \eqno(11\mbox{c})$$
The system (11c) admits many solutions. We here consider only the solution
$$B_1 = 2 h, \qquad B_4 = h, \qquad B_6 = 0 \qquad \mbox{and} \qquad
  A_{11} = 1. \eqno(11\mbox{d})$$
In this case, the relations (9) are of the following form
$$[u, {\sf d}u] = 2h {\sf d} u, \qquad [\eta, {\sf d}u]_+ = 0, $$
$$[u,{\sf d}\eta] = h {\sf d} \eta, \qquad [\eta, {\sf d}\eta] = 0. \eqno(12)$$
Applying the exterior differential {\sf d} to the first and second (or third)
of the relations (12) one gets
$$[{\sf d} u, {\sf d}\eta] = 0, \qquad ({\sf d} u)^2 = 0. \eqno(13)$$

{\it Note}. The two superalgebras (1) and (7) are closely related.
Therefore the differential calculi on these two superalgebras are
also closely related. Indeed, the differentials of $u$ and $\eta$
in terms of $x$ and $\theta$ are
$${\sf d} u = {{2h}\over {e^{2h} - 1}} ~{\sf d}x ~x^{-1}, \qquad
{\sf d} \eta = x^{-1} ({\sf d} \theta - {\sf d}x x^{-1} \theta)
\eqno(14)$$ so that replacing these into (12) one obtains
$$x ~{\sf d}x = q^2 {\sf d}x ~x, \qquad \theta ~{\sf d}x = - q {\sf
d}x ~\theta,$$
$$x ~{\sf d} \theta = q {\sf d} \theta ~x + (q^2 - 1) {\sf d}x
~\theta, \qquad \theta ~{\sf d} \theta = {\sf d} \theta ~\theta.$$
This differential structure is invariant under action of
$GL_q(1\vert 1)$ ~(see, e.g. [13]). Thus the differential
structure (12) must be invariant under action of $gl_h(1\vert 1)$
with $q = e^h$ (see, Ref. 11).

A differential algebra of the associative algebra ${\cal B}$ is a
Z$_2$-graded associative algebra $\Gamma$ equipped with an
operator {\sf d} that has the properties (10). Furthermore, the
algebra $\Gamma$ has to be generated by $\Gamma^0 \cup \Gamma^1
\cup \Gamma^2$, where $\Gamma^0$ is isomorphic to ${\cal B}$. For
${\cal B}$ we write ${\cal L(A)}$, the Lie superalgebra of $\cal
A$. Let us denote the algebra generated by ${\sf d} u$ and ${\sf
d} \eta$ with the relations (12) by $\Gamma^1$, where $\Gamma^1$
is isomorphic to ${\sf d} {\cal L(A)}$, and the algebra (13) by
$\Gamma^2$. Let $\Gamma$ be the quoitent algebra of the free
associative algebra on the set $\{u,\eta,{\sf d}u,{\sf d}\eta\}$
modulo the ideal $J$ that is generated by the relations (7), (12)
and (13). Then the differential algebra $\Gamma$ is a Z$_2$-graded
Hopf algebra with the following costructures:
$$\Delta({\sf d} u) = {\sf d} u \otimes 1 + 1 \otimes {\sf d} u, \qquad
  \Delta({\sf d} \eta) = {\sf d} \eta \otimes 1 + 1 \otimes {\sf d} \eta.$$
$$\epsilon({\sf d}u) = 0, \qquad \epsilon({\sf d}\eta) = 0,  \eqno(15)$$
$$\kappa({\sf d} u) = - {\sf d} u, \qquad
  \kappa({\sf d} \eta) = - {\sf d} \eta. $$

Before closing this section, just as we introduced the derivatives of the
generators of $\cal A$ in the standard way, let us introduce derivatives of
the generators of $\cal L(A)$ and multiply explicit expression of the exterior
differential {\sf d} from the right by $u f$ and $\eta f$, respectively. Then,
using the Z$_2$-graded Leibniz rule for partial derivatives
$$\partial_i (f g) = (\partial_i f) g + (-1)^{\hat{f}} f (\partial_i g)
  \eqno(16)$$
we get
$$[\partial_u, u] = {{2h}\over {e^{2h} - 1}} + 2 h \partial_u, \qquad
  [\partial_u, \eta] = 0,$$
$$[\partial_\eta, u] = h \partial_\eta, \qquad
  [\partial_\eta, \eta]_+ = 1. \eqno(17)$$
The commutation relations between the derivatives can be easily obtained by
using ${\sf d}^2 = 0$. So it follows that
$$0 = {\sf d}^2 = {\sf d} u ~{\sf d} \eta (\partial_u \partial_\eta -
  \partial_\eta \partial_u) + ({\sf d}\eta)^2 \partial_\eta^2$$
which says that
$$[\partial_u, \partial_\eta] = 0, \qquad \partial_\eta^2 = 0. \eqno(18)$$
Finally to find the commutation relations between the differentials and
derivatives we shall assume that they have the following form
$$\partial_u {\sf d}u = C_{11} {\sf d}u \partial_u +
   C_{12} {\sf d}\eta \partial_\eta + D_1 {\sf d}u + D_2 {\sf d} \eta, $$
$$\partial_u {\sf d} \eta = C_{21} {\sf d} \eta \partial_u +
   C_{22} {\sf d} u \partial_\eta + D_3 {\sf d}u + D_4 {\sf d} \eta, $$
$$\partial_\eta {\sf d}u = F_{11} {\sf d}u \partial_\eta +
   F_{12} {\sf d}\eta \partial_u + D_5 {\sf d}u + D_6 {\sf d} \eta, $$
$$\partial_\eta {\sf d} \eta = F_{21} {\sf d} \eta \partial_\eta +
   F_{22} {\sf d} u \partial_u + D_7 {\sf d}u + D_8 {\sf d} \eta. \eqno(19)$$
After some tedious but straightword calculations, we find
$$\partial_u {\sf d}u = e^{-2h} {\sf d}u \partial_u - e^{-2h} {\sf d}u,
  \qquad
 \partial_u {\sf d}\eta = e^{-2h} {\sf d}\eta \partial_u-e^{-2h} {\sf d}\eta,$$
$$\partial_\eta {\sf d} u = - {\sf d} u \partial_\eta, \qquad
  \partial_\eta {\sf d} \eta = {\sf d} \eta \partial_\eta +
  {{e^h - e^{-h}}\over {2h}} \left\{ (e^h - e^{-h}) {\sf d} u \partial_u +
  e^{-h} {\sf d} u \right\}. \eqno(20)$$

\noindent
{\bf 4. Hopf algebra structure of forms on $\cal L(A)$}

\noindent
Using the generators of $\cal L(A)$ we can define two one-forms as follows:
$$\phi = {{e^{2h} - 1}\over {2h}} {\sf d} u, \qquad V = e^h {\sf d} \eta.
  \eqno(21)$$
We denote the algebra of forms generated by two elements $\phi$ and $V$ by
$\Omega$. The generators of the algebra $\Omega$ with the generators of
$\cal L(A)$ satisfy the following relations:
$$[u, \phi] = 2h \phi, \qquad [\eta, \phi]_+ = 0, $$
$$[u, V] = h V, \qquad [\eta, V] = 0. \eqno(22)$$
The commutation rules of the generators of $\Omega$ are
$$\phi^2 = 0, \qquad [\phi, V] = 0. \eqno(23)$$

One can make the algebra $\Omega$ into a Z$_2$-graded Hopf algebra with the
following co-structures: the coproduct
$\Delta: \Omega \longrightarrow \Omega \otimes \Omega$ is defined by
$$\Delta(\phi) = \phi \otimes 1 + 1 \otimes \phi, \qquad
  \Delta(V) = V \otimes 1 + 1 \otimes V. \eqno(24) $$
The counit $\epsilon: \Omega \longrightarrow {\cal C}$ is given by
$$\epsilon(\phi) = 0, \qquad \epsilon(V) = 0 \eqno(25)$$
and the coinverse $\kappa : \Omega \longrightarrow \Omega$ is
defined by
$$\kappa(\phi) = - \phi, \qquad \kappa(V) = - V. \eqno(26)$$
One can easily to check that (22) and (23) are satisfied. Note
that the commutation relations (22) and (23) are compatible with
$\Delta$, $\epsilon$ and $\kappa$, in the sense that $\Delta(u
\phi) = \Delta(\phi u) + 2h \Delta(\phi)$ and so on.

\noindent
{\bf 5. The Hopf superalgebra of vector fields on $\cal
L(A)$}

\noindent
In this section, we shall obtain the superalgebra of vector fields
on  $\cal L(A)$ and their Hopf algebra structure. We first write
the Cartan-Maurer forms as
$${\sf d} u = {{2h}\over {e^{2h} - 1}} \phi, \qquad {\sf d} \eta = e^{-h} V.
  \eqno(27)$$
Then, the relations (23) allow us to calculate the superalgebra of
the vector fields. Writing the exterior differential {\sf d} in
the from
$${\sf d} = \phi X + V \nabla \eqno(28)$$
and considering an arbitrary function $f$ of the generators of $\cal A$ and
using the nilpotency of {\sf d} one has
$$\phi {\sf d} X = V {\sf d} \nabla. \eqno(29)$$
So we find the following commutation relations for the
superalgebra of vector fields
$$[X, \nabla] = 0, \qquad \nabla^2 = 0. \eqno(30)$$

We also note that the commutation relations (30) of the vector
fields should be consistent with monomials of the generators of
$\cal L(A)$. To proceed, we must calculate the actions of the
Z$_2$-graded Leibniz rule by comparing the elements which lie
together with each other from the one-forms:
$$[X, u] = 1 + 2h X, \qquad [X, \eta] = 0, $$
$$[\nabla, u] = h \nabla, \qquad [\nabla, \eta]_+ = 1.
  \eqno(31)$$
Of course, these commutation relations must be consistent.

In order to find the coproduct of this superalgebra, we shall use
the fact that the exterior differential operator {\sf d} satisfies
the Leibniz rule [13]. So using the Z$_2$-graded Leibniz rule
for {\sf d} we write
$$X(f g) = (X f) g + e^{2hN} f (Xg), $$
$$\nabla(f g) = (\nabla f) g + e^{hN} f (\nabla g), $$
where $N$ is a number operator which acting on the monomials of
the generators of ${\cal L(A)}$. This provides a comultiplication
$$\Delta(X) = X \otimes I + e^{2hN} \otimes X, $$
$$\Delta(\nabla) = \nabla \otimes I + e^{hN} \otimes \nabla. \eqno(32)$$

Using the following basic axioms of Hopf superalgebra
$$m(\epsilon \otimes \mbox{id}) \Delta(Y) = Y, \qquad
  m(\mbox{id} \otimes \kappa) \Delta(Y) = \epsilon(Y) \eqno(33)$$
one obtains
$$\epsilon(X) = 0, \qquad \epsilon(\nabla) = 0, \eqno(34)$$
$$\kappa(X) = - e^{-2hN} X, \qquad
  \kappa(\nabla) = - e^{-hN} \nabla. \eqno(35)$$

We can now easily obtain the dual Hopf superalgebra as follows: if
we introduce the operators $N$ and $\chi$ as
$$e^{2hN} = I + (e^{2h} - 1) X, \qquad
  \chi = \nabla \left\{I + (e^{2h} - 1) X\right\}^{-1/2} \eqno(36)$$
then we have
$$\Delta(N) = N \otimes I + I \otimes N, $$
$$\Delta(\chi) = \chi \otimes e^{-hN} + I \otimes \chi. \eqno(37\mbox{a})$$
and
$$\epsilon(N) = 0, \qquad \epsilon(\chi) = 0, \eqno(37\mbox{b})$$
$$\kappa(N) = - N, \qquad \kappa(\chi) = - \chi e^{hN}. \eqno(37\mbox{c})$$
The operators $N$ and $\chi$ are preserve the commutation
relations (30):
$$[N, \chi] = 0, \qquad \chi^2 = 0. \eqno(38)$$

Note that the relations (38) and the Hopf algebra structure (37)
can also be obtained from the approach of Ref. 11 [see, chapter
3].

An interesting problem is the construction of a differential
calculus on the Lie superalgebra of the quantum supergroup GL$_q(1\vert 1)$
using the methods of this paper and Ref. 14.

{\bf Acknowledgement}

This work was supported in part by TBTAK the turkish Scientific
and Technical Researc Council.

\vfill\eject 

\end{document}